\documentclass{amsart}

\usepackage{geometry}
\usepackage{color, graphicx}
\usepackage[all]{xy}

\usepackage{calrsfs}
\usepackage{times}
\usepackage[scaled=0.92]{helvet}

\usepackage{amsthm}
\theoremstyle{plain}

\theoremstyle{definition}

\swapnumbers
\theoremstyle{plain}
\newtheorem{thm}{Theorem}[section]
\newtheorem{prop}[thm]{Proposition}
\newtheorem{cor}[thm]{Corollary}
\newtheorem{lem}[thm]{Lemma}
\newtheorem{prob}[thm]{Subgroup Conjugacy Problem}
\theoremstyle{definition}
\newtheorem{df}[thm]{Definition}

\newtheorem{rem}[thm]{Remark}
\newtheorem{se}[thm]{}
\newtheorem*{ack}{Acknowledgments}

\newcommand{\abs}[1]{\lvert#1\rvert}
\def\ZR#1{\langle #1\rangle}

\usepackage{amssymb}
\usepackage{bbm}

\def\Z{\mathbf{Z}}

\def\R{\mathbf{R}}

\def\F{\mathbf{F}}
\def\P{\mathbf{P}}

\def\Aut{\mathop{\mathrm{Aut}}\nolimits}

\def\Gal{\mathop{\mathrm{Gal}}\nolimits}

\def\PGL{\mathop{\mathrm{PGL}}\nolimits}

\def\PSL{\mathop{\mathrm{PSL}}\nolimits}

\begin{document}

\date{\today\ (version 1.0)} 
\title[Relating deformation parameters]{Three examples of the relation between rigid-analytic\\ and algebraic deformation parameters}

\author[G.~Cornelissen]{Gunther Cornelissen}
\address{Mathematisch Instituut, Universiteit Utrecht, Postbus 80.010, 3508 TA Utrecht, Nederland}
\email{g.cornelissen@uu.nl}
\author[F.~Kato]{Fumiharu Kato}
\address{Department of Mathematics,  Faculty of Sciences, University of Kyoto, Kyoto 606-8502, Japan}
\email{kato@math.kyoto-u.ac.jp}
\author[A.~Kontogeorgis]{Aristides Kontogeorgis}
\address{Department of Mathematics, University of the Aegean, Karlovasi, Samos, Greece}
\email{kontogar@aegean.gr}

\begin{abstract} 
\noindent We consider three examples of families of curves over a non-archimedean valued field which admit a non-trivial group action. These equivariant deformation spaces can be described by algebraic parameters (in the equation of the curve), or by rigid-analytic parameters (in the Schottky group of the curve). We study the relation between these parameters as rigid-analytic self-maps of the disk. \end{abstract}

\maketitle

\section*{Introduction} We consider a situation where one has two natural sets of deformation parameters for the same deformation problem, and ask what the relation between these two sets of parameters is. In our case, this leads to interesting analytic relations, generalising the relation between the Legendrian and the Tate parameter on an elliptic curve. 

More precisely, the situation we study is as follows. Suppose $X \rightarrow Y$ is a Galois cover of projective curves over a non-archimedean valued field $K$. Suppose $X$ is a Mumford curve over $K$ with Schottky group $\Gamma$. Then the situation can be described by a diagram
$$ \xymatrix{   \Omega_\Gamma \ar[ddr]_N \ar[dr]^\Gamma & \\ & X \ar[d]^{G=\Gamma \backslash N} \\ & Y}$$ 
where $\Omega_\Gamma$ is the domain of discontinuity of $\Gamma$ in $\P^1(K)$ and $N$ is a subgroup of $\PGL(2,K)$ containing $\Gamma$ and contained in the normalizer of $\Gamma$. Thus, $N$ gives a kind of ``orbifold'' uniformization of $Y$. 

Now suppose the situation in this diagram admits moduli: in the algebraic sense it means that $X$ with the action of $G$ has nontrivial deformations; in the analytic sense it means that the embedding of $N$ in $\PGL(2,K)$ admits nontrivial deformations (cf.\ e.g.\ \cite{CoKa2} for formal definitions). We know from \cite{CoKa2}, \S 9, that the formal algebraic and analytic deformation functors (defined in that reference) are isomorphic. But what is the nature of the relation between the analytic and algebraic deformation parameters? It can hardly be algebraic in the general case. We give an explicit form of this map as a rigid analytic infinite product in three cases. For Tate elliptic curves $X$ covering $Y=\P^1$, we find back a classical relation in terms of sigma functions (cf.\ \ref{tate}). For four point covers of $\P^1$ by Mumford curves, we find a `localisation formula' for the cross ratio (cf.\ \ref{four}). The final, most elaborate, case is over a field $K$ of characteristic $p>0$. Suppose $K$ contains the finite field $\F_q$ with $q=p^m$. Let $\Delta_0$ denote the punctured open unit disk $\Delta_0 = \{ 0 < |z| < 1 \}$, and $\Delta$ the open unit disk $\Delta = \{ |z|<1 \}$. For $\lambda \in K$, let $X_{\lambda}$ denote the algebraic curve 
$$ X_{\lambda} = \{ (x^q-x)(y^q-y)=\lambda \} \subseteq \P^1 \times \P^1. $$
For $\lambda \neq 0$, this is a smooth curve of genus $g:=(q-1)^2$, and for $\lambda \in \Delta_0$, it defines a Mumford curve. 

These curves are quite special from a few points of view. They are the Mumford curves with maximal automorphism group (and hence their Schottky groups are the analogue of classical Hurwitz groups), cf.\ \cite{CKK}, \cite{CoKa3}. They are highly symmetric ordinary curves and were first studied by  Subrao (\cite{Subrao}), Valentini-Madan (\cite{Val-Mad:80} and Nakajima (\cite{Nak}) from this point of view.  Michel Matignon has studied their equivariant liftability to characteristic zero \cite{Mat}, and the curves play a special role when studying the field of definition versus field of moduli question for cyclic covers of the projective line (cf. \cite{Kon}).

A corresponding Schottky group can be constructed as follows: define elements of $\PGL(2,K)$  by $$
\varepsilon_u:=\bigg(\begin{array}{cc}1&u\\ 0&1\end{array}\bigg),\quad\tau:=\bigg(\begin{array}{cc}0&t\\ 1&0\end{array}\bigg),\quad \varepsilon'_u:=\tau\varepsilon_u\tau
$$
where $u\in\F_q$ and $t\in K^{\times}$ with $t \in \Delta_0$. Note that these are explicit matrices whose entries are polynomial in $t$ with $\F_q$-coefficients. If we define a subgroup $\Gamma=\Gamma(t)$ of $\PGL(2,K)$ to be generated by the commutators of these elements:
$$\Gamma(t):=\ZR{[\varepsilon_u,\varepsilon'_v]\,|\,u,v\in\F_q},$$ then it is known that the Schottky group of $X_\lambda$ occurs as $\Gamma(t)$ for some $t \in \Delta_0$.
 
The algebraic cover $X_\lambda \rightarrow \P^1$ is ramified over a fixed set of 3 points. Nevertheless, it admits non-trivial equivariant moduli (this is a typical phenomenon in positive characteristic and cannot occur in characteristic zero). On the algebraic side these moduli are given by varying $\lambda$. In some sense, this $\lambda$ is a substitute for the cross ratio if the number of branch points is below four. We prove that $X_\lambda$ and $X_\mu$ are isomorphic precisely if $\lambda/\mu \in \F_q^*$ (Theorem \ref{subraocurves}). 

On the analytic side, the moduli are given by varying $t$. We prove that two Schottky groups $\Gamma(t)$ and $\Gamma(u)$ are conjugate precisely if $t/u \in \F_q^*$ (Theorem \ref{asmconj}). We observe that the method of proof also provides an amusing `analytic' solution to the conjugacy problem for subgroups of $\PGL(2,K)$ generated by matrices $\varepsilon_u, \varepsilon'_v$ as above (in general, such problems are considered quite hard). 

Our further main theorems (cf.\ \ref{mainthm}, \ref{analytic},  \ref{mainthm2}) will provide an expression for $\lambda$ as an infinite product of rational functions of $t$, that defines a rigid-analytic self-map on the punctured disk $\Delta_0$ and hence extends to $\Delta$. We also give a rigid analytic \emph{automorphism} of $\Delta_0$ that provides a correspondence between the actual moduli, so between isomorphism classes of $X_\lambda$ and conjugacy classes of $\Gamma(t)$.

\begin{ack}
We thank Andreas Schweizer for help with the conjugacy problem for ASM-groups. 
\end{ack}

\section{Infinite products}

Before we get started, we recall the following lemma from $p$-adic analysis. Let $K$ denote a non-archimedean valued field. An open set of $\P^1$ is seen as an analytic subspace of $\P^1$ as in \cite{FvdP}, Def.\ 2.6.1 (I.7 in the original French edition). Also, ``holomorphic function'' should be understood as in that reference, section 2.2. 
\begin{lem} \label{analyticlem}
\begin{enumerate}
\item
A product $\prod\limits_{n=1}^\infty a_n$ converges if and only if $\lim\limits_{n \rightarrow \infty} a_n =1$.
\item
An infinite product  $\prod\limits_{n=1}^\infty a_n(t)$ of holomorphic functions $a_n(t)$
converges uniformly on an open set $U \subset \P^1(K)$ to an analytic function, if there is a sequence $c_n$ of real numbers such that 
$\left| a_n(t) -1\right| \leq c_n$ for all $t \in U$ with  $\lim \limits_{n \rightarrow \infty} c_n=0$. 
\item If $f$ is a bounded holomorphic function on the punctured open disk $\Delta_0$ then it extends to a holomorphic function 
on the whole open disk $\Delta$.
 \end{enumerate}
\end{lem}
\begin{proof}
The first two statements are analogues of classical results in complex analysis, see, e.g.\ \cite{SS} ch.\ 5, section 3, but easier since in non-archimedean analysis, a sum converges if and only if the general term tends to zero. The third statement follows from Remark 2.7.14 in \cite{FvdP} (or from the Lemma after I.8.11 in the original French edition), applied to closed disks $\{ 0<|t|<1-\epsilon \}$ before taking $\epsilon \rightarrow 0$.
\end{proof}

\section{Tate elliptic curve}

\begin{se} In this section, we present a well understood situation as a guiding example for what is to come in later paragraphs. 

Consider a Tate curve $E_q = {\bf G}_m(K)/q^\Z$ over a non-archimedean valued field $K$. It covers $\P^1$ over four points. Without loss of generality, we can assume the branch points are $\{0,1,\infty,\lambda\}$ with $|\lambda|=1$. It is known that $\lambda$ then satisfies $|1-\lambda|<|2|^2$ (cf.\ Bradley, \cite{Bradley}, Example 3.8). We now compute the (well-understood) relationship between the period $q$ and the Legendrian 2-torsion modulus $\lambda$ as a toy model for the method we will use later on. 

We let $N$ denote the group generated by $\alpha:=\bigl( \begin{smallmatrix} 0 & 1 \\ 1 & 0 \end{smallmatrix} \bigr)$ and $\beta:=\bigl( \begin{smallmatrix} 0 & t^2 \\ 1 & 0 \end{smallmatrix} \bigr)$ and we let $\Gamma$ be the group generated by $\alpha \beta = \bigl( \begin{smallmatrix} 1 & 0 \\ 0 & t^2 \end{smallmatrix} \bigr)$. We set $q=t^{-2}$. Thus, $t$ corresponds to the choice of a two-torsion point in the Tate model $E_q$ belonging to the group $\Gamma$. The orbifold belonging to $N$ is the quotient of $E_q$ by an element of order two, leading to $\P^1$ with four branch points. 
\end{se} 
\begin{prop}\label{tate} The cross ratio of the four points of order two on the Tate curve $E_q$ with $q=t^{-2}$, for $|t|>1$ is given by
$$ \lambda(t) =\left(  \prod_{i \geq 0} \frac{t^{2i+1}+1}{t^{2i+1}-1} \right)^8. $$ 
The function $\lambda$ is a holomorphic function on the open disk $|t|>1$ centered at $\infty \in \P^1(K)$.
\end{prop}

\begin{rem} 
This is a classical formula in the theory of the Weierstrass $\sigma$-function, see \cite{AECII}, p. 89 (1.18).
\end{rem}

\begin{proof}
The $\Gamma$-orbit of the fixed points of $\alpha$ are $\{q^i\}$ and $\{-q^i\}$ for $i \in \Z$, which we suppose map to $0$ and $\infty$, respectively, on $\P^1$ (this means $|q|<1$, so $|t|>1$). Now the $N$-orbits of the fixed points are the same, with stabilizers of order two: the quotient $N/\Gamma$ is a group of order two generated by the class of $\alpha$, and its non-trivial element acts on $q^i$ by $q^i \mapsto \alpha \cdot q^i =q^{-i}$; on the algebraic side, this is the action $E_q \mapsto E_q / \{ \pm 1 \} = \P^1$, see the diagram in the introduction with $X=E_q$, $Y=\P^1$ and $G=\{\pm 1\}$. Therefore, if we define $u$ by
$$ u: = \kappa_1 \cdot \prod_{i \in \Z} \left( \frac{z-q^i}{z+q^i} \right)^2 $$
then $u$ is a parameter on the quotient $\P^1$. If we want $u$ to take on the value $1$ at $z=t$, this fixes the constant $\kappa_1$, so
$$ u(z) = \prod_{i \in \Z} \left( \frac{(t+q^i)(z-q^i)}{(t-q^i)(z+q^i)}\right)^2. $$
The fixed points of $\beta$ are $\{tq^i\}$ and $\{-tq^i\}$ for $i \in \Z$, which we suppose are above $1$ (at $z=t$) and $\lambda$ (at $z=-t$) on $\P^1$. Therefore, in the same way we find
$$ \lambda = \prod_{i \in \Z} \left( \frac{(t+q^i)(-t-q^i)}{(t-q^i)(-t+q^i)} \right)^2 =   \prod_{i \in \Z} \left( \frac{t^{2i+1}+1}{t^{2i+1}-1} \right)^4. $$
(Note that the product here is over $i \in \Z$.)

For proving that $\lambda(t)$ is a holomorphic function of $t$ we are going to use Lemma \ref{analyticlem}.
Observe that for $|t|>c>1$,
\begin{equation} \label{l-1}
\left| \frac{t^{2i+1}+1 }{ t^{2i+1}-1 }-1\right|=
\frac{\left|2\right|}{\left| t^{2i+1}-1\right|}=\frac{|2|}{|t|^{2i+1}} \leq \frac{1}{c^{2i+1}} \stackrel{i \rightarrow \infty}{\longrightarrow} 0.
\end{equation}
Therefore, the product defining $\lambda(t)$ is uniformly convergent on $|t|>c>1$. We are hence allowed to reorder the terms in the product expression for $\lambda(t)$ from taking the product over $i \in \Z$ to $i \geq 0$ and squaring, to arrive at the formula indicated in the proposition. 
\end{proof}

\begin{se} We now compare isomorphism of elliptic curves expressed in $t$ and $\lambda$. 
Two Tate elliptic curves $E_{t_1},E_{t_2}$ 
are isomorphic if the free groups 
$\Gamma_{t_i}=\langle \bigl( \begin{smallmatrix} 1 & 0 \\ 0 & t_i^2 \end{smallmatrix} \bigr) \rangle$  uniformizing them
are conjugate. This means their (two possible) generators should map to each other. By direct computation, one finds that this is equivalent to $t_1=\pm t_2$. So the moduli of Tate curves in terms of the $t$-parameter is given by dividing the disk $|t|>1$ by the automorphism $t \mapsto -t$. 

From the above infinite product expansion, we see that everything is indeed normalized so that $|\lambda(t)|=1$. We also find immediately the behaviour of $\lambda$ w.r.t.\ the disk automorphism $t \mapsto -t$, namely $\lambda(-t)=\lambda(t)^{-1}$. Also observe that the Mumford conditions $|\lambda|=1$ and $|1-\lambda|<|2|^2$ are stable under $\lambda \mapsto \lambda^{-1}$.

On the algebraic side, two elliptic curves are isomorphic if they have the same $j$-invariant.
We can express the $j$-invariant as a function of the cross ratio $\lambda$ by 
\[
j(\lambda)=2^8 \frac{(\lambda^2-\lambda+1)^3}{\lambda^2(\lambda-1)^2}.
\]
(The Mumford conditions arise from putting $|\lambda|=1$ in this expression and requiring $|j(\lambda)|>1$.) This defines the $S_3$-cover of modular curves $Y_0(2) \rightarrow Y(1)$.
For every $j$-invariant there are generically  $6$ values of $\lambda$ mapping to this $j$-invariant, 
namely $$\lambda,\frac{1}{\lambda},1-\lambda,\frac{1}{\lambda-1},\frac{\lambda}{\lambda-1},
\frac{\lambda-1}{\lambda},$$ but the only one of these that satisfies the Mumford conditions (given that $\lambda$ satisfies them) is $\lambda^{-1}$, as is seen by direct verification. 

The infinite product $\lambda(t)$ induces a bijection 
\begin{eqnarray*} \{ |t|>1 \} /\{ t \mapsto -t \}  & \rightarrow & \{ |\lambda|=1 \mbox{ and } |1-\lambda|<1 \}/\{ \lambda \mapsto \lambda^{-1} \}  \\  T & \mapsto &  \Lambda(T):=\lambda^{\pm 1} (\pm T), \end{eqnarray*}
which is now one-to-one on isomorphism classes of Tate curves. 

Compare this section with the discussion in \cite{BGR}, 9.7.3. 
\end{se}

\section{Four point covers of $\P^1$ by Mumford curves}
In this section we show a relation between global and local cross ratios on Mumford curves. Again, this is a kind of toy model for the theory.
Let $\pi: X \rightarrow \P^1$ be a covering of Mumford curves that is branched above four points $0,1,\infty, \lambda$, let $\Gamma$ denote the Schottky group of $X$ and $N$ the subgroup of $\PGL(2,K)$ corresponding to the covering $\pi$. Let $E_x$ denote the set of fixed points on $\P^1(K)$ of the elliptic elements in $N$ that are above $x \in \{0,1,\infty,\lambda\}$ and let $e_x \in E_x$ denote one chosen fixed element for each such $x$. For $\kappa_1 \in K^*$, the function 
$$ u(z) := \kappa_1 \cdot \prod_{\gamma \in N} \frac{z-\gamma(e_0)}{z-\gamma(e_\infty)} $$
is a meromorphic function of $z$ with simple zeros at $E_0$ and simple poles at $E_\infty$ (cf.\ \cite{GvdP}, pp. 44--47). Therefore, it is a uniformizer on $\P^1$. We normalize it so it takes the value $1$ at $e_1$; this determines the value of $$\kappa_1 = \prod \frac{e_1-\gamma(e_\infty)}{e_1-\gamma(e_0)} $$ uniquely, and then $u(e_\lambda)=\lambda$ gives the following expression for $\lambda$:
$$ \lambda = \prod_{\gamma \in N} \frac{e_\lambda-\gamma(e_0)}{e_\lambda-\gamma(e_\infty)} \cdot \frac{e_1-\gamma(e_\infty)}{e_1-\gamma(e_0)}.$$
Since on the right hand side the factors in the product are cross ratios, we can rewrite this formula as follows:

\begin{prop} \label{four} With these notations, 
$$ \lambda = (0,1;\infty,\lambda) = \prod_{\gamma \in N} (e_\lambda, e_1; \gamma(e_0), \gamma(e_\infty)).$$
\end{prop}

\begin{rem}
Whether a four-point cover of $\P^1$ is a Mumford curve depends on the location of the branch points (so on $\lambda$) and is a rather subtle question. For some answers, see Bradley \cite{Bradley}, Theorem 5.4. 
\end{rem}

\section{Artin-Schreier-Mumford group}\label{sec-ASMcurves}
In the rest of this paper, we analyse our main example: certain fiber products of $p$-covers in characteristic $p$ that admit non-trivial equivariant moduli, but with fixed branch points. Hence there is no technology of cross ratios to rely upon. As a matter of fact, this work could be seen as constructing a kind of cross ratio for wild covers with less than four branch points. We work in the following situation:
\begin{itemize}
\item $K$ is a complete non-archimedean valued field of characteristic $p>0$;
\item $\abs{\,\cdot\,}\colon K\rightarrow\R_{\geq 0}$ is a multiplicative valuation (that is, a norm) of $K$;
\item $q=p^m$, where $m>0$ is a positive integer;
\item we assume that $K$ contains $\F_q$.
\end{itemize}

We often write fractional linear transformations, that is, elements of $\PGL(2,K)$, in matrix form.
Set
$$
\varepsilon_u:=\bigg(\begin{array}{cc}1&u\\ 0&1\end{array}\bigg),\quad\tau:=\bigg(\begin{array}{cc}0&t\\ 1&0\end{array}\bigg),
$$
where $u\in\F_q$ and $t\in K^{\times}$.
Notice that $\tau$ is of order $2$, and that the map $u\mapsto\varepsilon_u$ is an injective group homomorphism from the additive group $\F_q$ into $\PGL(2,K)$; in particular, we have $\varepsilon_u\varepsilon_v=\varepsilon_{u+v}$ for $u,v\in\F_q$ and $\varepsilon^{-1}_u=\varepsilon_{-u}$.
We set 
$$
\varepsilon'_u:=\tau\varepsilon_u\tau=\bigg(\begin{array}{cc}t&0\\ u&t\end{array}\bigg)
$$
for $u\in\F_q$.
The subgroup of $\PGL(2,K)$ consisting of all elements of the form $\varepsilon_u$ (resp.\ $\varepsilon'_u$) for $u\in\F_q$ is denoted by $E$ (resp.\ $E'$).
The subgroups $E$ and $E'$ are elementary abelian subgroups, isomorphic to $(\Z/p\Z)^m$.

For any $a\in K^{\times}$, we define 
$$
\mu_a:=\bigg(\begin{array}{cc}a&0\\ 0&1\end{array}\bigg), 
$$
that is, $\mu_a(z)=az$.
The following identities are easy to see:
$$
\varepsilon_u\mu_v=\mu_v\varepsilon_{u/v},\quad \varepsilon'_u\mu_v=\mu_v\varepsilon'_{uv},\quad \tau\mu_v=\mu_{v^{-1}}\tau,
$$
for $u\in\F_q$ and $v\in\F^{\times}_q$.

We consider the following subgroups of $\PGL(2,K)$
\begin{itemize}
\item $N_1:=\ZR{\varepsilon_u,\varepsilon'_v\,|\,u,v\in\F_q}$;
\item $\Gamma:=\ZR{\delta_{u,v}\,|\,u,v\in\F_q}$,
\end{itemize}
where $\delta_{u,v}$ denotes the commutator of $\varepsilon_u$ and $\varepsilon'_v$, that is, 
$$
\delta_{u,v}:=[\varepsilon_u,\varepsilon'_v]=\varepsilon_u\varepsilon'_v\varepsilon_{-u}\varepsilon'_{-v}.
$$
The method of isometric circles (for example) proves that for $\abs{t}<1$, these groups are discrete in $\PGL(2,K)$ (cf.\ \cite{CKK}, \S 8). 

Notice that the following identity holds:
$$
\delta^{-1}_{u,v}=\tau\delta_{v,u}\tau.
$$

Clearly, $\Gamma$ is a normal subgroup of $N_1$ such that
$$
N_1/\Gamma\cong E\times E'\cong(\Z/p\Z)^{2m}.
$$
It is known (cf.\ \cite{CKK}, \S 9) that the subgroup $N_1$ is isomorphic to the free product $E\ast E'$; moreover, the subgroup $\Gamma$ is a free group of rank $(q-1)^2$, and has exactly $\delta_{u,v}$ for $u,v\in\F^{\times}_q$ as its free generators.

\begin{rem}\label{rem-normalizer}{\rm 
The normalizer $N=N(\Gamma)$ of $\Gamma$ in $\PGL(2,K)$ is generated by $E$, $\tau$, and $\mu_v$ for $v\in\F^{\times}_q$.}
\end{rem}

We set
$$
N_0:=\ZR{\Gamma,E'},\quad N'_0:=\ZR{\Gamma,E},
$$
which are subgroups of $N_1$ containing $\Gamma$.
\begin{lem}\label{lem-section}
We have $N_0/\Gamma\cong E'$ and $N'_0/\Gamma\cong E$.
\end{lem}

\begin{proof}
The assertion is clear, since the first and second factors of $N_1/\Gamma\cong E\times E'$ have sections $E',E\hookrightarrow N_1$.
\end{proof}

We end this section with a diagram of groups:

\begin{equation} \label{dia} \xymatrix{0 \ar[r] & \Gamma \ar[r] & N \ar[r] & (E \times E') \rtimes D_{q-1} \ar[r] & 0  \\
0 \ar[r] & \Gamma \ar[r] \ar@{=}[u] &  N_1  \ar@{^{(}->}[u] \ar[r] & E \times E' \ar[r] \ar@{^{(}->}[u] & 0  \\
0 \ar[r] & [E,E']  \ar@{=}[u] \ar[r]  & E \ast E'  \ar[r] \ar@{=}[u] &  E \times E' \ar@{=}[u] \ar[r] & 0
}\end{equation} 

\section{Artin-Schreier-Mumford curves}\label{sec-ASMcurves1}
Let $\Omega\subset\P^1_K$ be the set of ordinary points with respect to the Schottky group $\Gamma$.
The quotient $X:=\Omega/\Gamma$ is a Mumford curve over $K$, which we call an {\em Artin-Schreier-Mumford curve (ASM-curve or ASM-cover)}.
Let us denote the quotient map $\Omega\rightarrow X$ by $\varphi$.
Similarly, one can consider the quotients $\Omega/N_0$, $\Omega/N'_0$, and $\Omega/N_1$, which we denote respectively by $X_0$, $X'_0$, and $X_1$.
We have the maps $X\rightarrow X_0$ and $X\rightarrow X'_0$, which are Galois covering maps with the Galois groups $E'$ and $E$, respectively; we denote these maps  by $\phi$ and $\phi'$, respectively.
Thus we get the following commutative diagram
$$
\xymatrix@R-3.5ex@C-2.2ex{&&X_0\ar[dr]\\ \Omega\ar[r]^{\varphi}\ar@/^1pc/[urr]^{\Phi}\ar@/_1pc/[drr]_{\Phi'}&X\ar[ur]^{\phi}\ar[dr]_{\phi'}&&X_1\\ &&X'_0\ar[ur]}
$$
(where we set $\Phi=\phi\circ\varphi$ and $\Phi'=\phi'\circ\varphi$).
By the proof of \cite[Prop.\ 2]{CoKa3}, we have:
\begin{itemize}
\item the curves $X_0$, $X'_0$, and $X_1$ are isomorphic to $\P^1_K$;
\item the square in the above diagram is Cartesian, that is, $X\cong X_0\times_{X_1}X'_0$.
\end{itemize}

Let $\zeta$ be an inhomogeneous coordinate of $X_1$.
Then, by [loc.\ cit.], we can choose homogeneous coordinates $(X:Z)$ for $X_0$, and $(Y:W)$ for $X'_0$, with corresponding inhomogeneous coordinates  $x=X/Z$ and $y=Y/W$, in such a way that \begin{itemize} \item the Galois covering map $X_0\rightarrow X_1$  is described by $x^q-x=\lambda\zeta$ for some $\lambda\in K^{\times}$ with $\abs{\lambda}<1$; \item the Galois covering map $X'_0\rightarrow X_1$ is given by $y^q-y=\zeta^{-1}$. \end{itemize}  
Hence the curve $X$, regarded as a closed subscheme in $X_0\times_KX'_0=\P^1_K\times_K\P^1_K$, is given by the explicit equation
\begin{equation} \label{*}
(X^q-XZ^{q-1})(Y^q-YW^{q-1})=\lambda Z^qW^q.
\end{equation}
Notice that $\lambda$ is uniquely determined, since the inhomogeneous coordinate $\zeta$ of $X_1$ was chosen so that the Galois covering map $X\rightarrow X_1$ branches over $\zeta=0,\infty$; another such choice differs only by a non-zero constant of $K$, and gives rise to the same $\lambda$.

In the sequel, we denote the curve $X$  by $X_{\lambda}$ (as a subscheme of $\P^1 \times \P^1$) to emphasize that it is given by the equation (\ref{*}). Note that since $\lambda \neq 0,1$, there is no confusion with $X_0,X_1$. 
By uniqueness of $\lambda$ as a parameter for this moduli problem, one can regard $\lambda$ as a function of $t$: $\lambda=\lambda(t)$. (Later on we will consider the effect of weakening the equivalence relation on $\lambda$ to that of mere isomorphism of abstract curves instead of embedded curves.) 

\section{Ramification points for ASM-covers} 
We know that the maps $\phi$ and $\phi'$ are Galois covering maps with Galois group isomorphic to $(\Z/p\Z)^m$; moreover, this action of $(\Z/p\Z)^m$ on $X_{\lambda}$ is expressed in the double homogeneous coordinate $(X:Z)\times(Y:W)$ as follows:
$$
\begin{array}{lcl}
\textrm{ for } \phi & : & (X:Z)\times(Y:W)\stackrel{u}{\longmapsto}(X:Z)\times(Y+uW:W)\\
\textrm{ for } \phi' & :  & (X:Z)\times(Y:W)\stackrel{u}{\longmapsto}(X+uZ:Z)\times(Y:W)
\end{array}
$$
for $u\in\F_q$.
By this, one easily see the following:
\begin{prop}\label{prop-ramifications1}
\begin{enumerate} \item The ramification points of the map $\phi$ are exactly the points of the form $(u:1)\times(1:0)$  for $u\in\F_q$. \item The ramification points of the map $\phi'$ are exactly the points of the form $(1:0)\times(u:1)$.
\item At each of these points, the map $\phi$ $($resp.\ $\phi')$ ramifies completely. \qed \end{enumerate}
\end{prop}

Next, to describe the maps $\Phi$ and $\Phi'$, let $z$ be the inhomogeneous coordinate of $\Omega\subset\P^1_K$ on which the  groups $\Gamma$, etc.\ act by fractional transformations in the given form.
\begin{prop}\label{prop-ramifications2}
\begin{enumerate} \item The ramification points of the map $\Phi$ are exactly the points of the form $\gamma(u)$ for $\gamma\in\Gamma$ and $u\in\F_q$. \item The ramification points of the map $\Phi'$ are exactly the points of the form $\gamma(t/u)$ for $\gamma\in\Gamma$ and $u\in\F_q$. 
\item At each of these points, the maps $\Phi$ $($resp.\ $\Phi')$ ramifies completely.
\item Moreover, by a change of coordinates $x \mapsto x+u$, $y \mapsto y+v$ $($for some  $u,v\in\F_q)$, one can assume that  for $u\in\F_q$, $\Phi$ maps the points $\gamma(u)$ to $(u:1)\times(1:0)$ and $\Phi'$ maps the points $\gamma(t/u)$ to $(1:0)\times(u:1)$.\end{enumerate}
\end{prop}

\begin{proof}
First, observe the following:
\begin{itemize}
\item for any $u,v\in\F_q$ $(v\neq 0)$, the element $\varepsilon_u\varepsilon'_v\varepsilon^{-1}_u=[\varepsilon_u,\varepsilon'_v]\varepsilon'_v$ belongs to $N_0$, being of order $p$, and has $z=u$ as its unique fixed point;
\item since $N_1$ is isomorphic to the free product $E \ast E'$, the normalizer of any non-trivial subgroup of $E'$ is $E'$ itself;
\item on the other hand, since $N_1/\Gamma\cong E\times E'$, the intersection $N_0\cap E'$ lies in $\Gamma$; but, since $\Gamma$ is a free group, we deduce that this intersection is $\{1\}$;
\item hence none of the elements of the form $\varepsilon_u\varepsilon'_v\varepsilon^{-1}_u$ are conjugate to each other.
\end{itemize}
Since the map $\Phi$ has to ramify exactly at $q$ points (that is, the points of the form $(u:1)$ for $u\in\F_q$) on $X_0$ and is totally ramified of order $q$ at each ramification point, elements of the form 
\begin{equation} \label{**}
\gamma\cdot\varepsilon_u\varepsilon'_v\varepsilon^{-1}_u\cdot\gamma^{-1}
\end{equation}
for $\gamma\in N_0$ and $u,v\in\F_q$ $(v\neq 0)$ are exactly the elements of finite order in $N_0$.
Since $\Gamma$ is a normal subgroup in $N_0$ and since one calculates
$$
\varepsilon'_w\cdot\varepsilon_u\varepsilon'_v\varepsilon^{-1}_u\cdot\varepsilon'_{-w}=[\varepsilon'_w,\epsilon_u]\cdot\varepsilon_u\varepsilon'_v\varepsilon^{-1}_u\cdot[\varepsilon'_w,\epsilon_u]^{-1},
$$
we may actually assume that the $\gamma$ in the expression (\ref{**}) are taken from $\Gamma$, and the first assertion follows (the assertion for $\Phi'$ is similar).

To show the last assertion, we only need to show that the map $\Phi$ is equivariant w.r.t.\ the respective actions of  $E$ on $\Omega$ and of $N_1/N_0\cong E$ on $X_0$, which is clear due to existence of the section $E\hookrightarrow N_1$ as in the proof of Lemma \ref{lem-section}.
\end{proof}

\begin{rem}\label{rem-NBordinary}{\rm 
By Proposition \ref{prop-ramifications2}, in particular, we know that the points of the form $\gamma(u)$ or $\gamma(t/u)$ (for $\gamma\in\Gamma$ and $u\in\F_q$) are in $\Omega$, and hence, are not fixed points of any non-trivial element of $\Gamma$.}
\end{rem}

\section{Evaluation of $\lambda(t)$}\label{sec-evaluation}

After the identification of ramification points in the analytic and the algebraic representation of the ASM-curve, we now find an analytic relation between $t$ and $\lambda$ roughly as follows: we are going to write the algebraic coordinate functions $x$ and $y$ as analytic functions of $t$ (using the location of their zeros and poles, we can do this with an infinite product). We then insert these expressions into the algebraic equation and analyse what happens analytically. 

\begin{se}[Coordinate functions]
We are in the following situation:
$$
\xymatrix@C-2ex{\Omega\ar[d]_{\varphi}\\ X_{\lambda}\ar@{^{(}->}[r]&\P^1_K\times_K\P^1_K,}
$$
where $\varphi$ is the quotient map by the free group $\Gamma$, and $\P^1_K\times_K\P^1_K$ is considered with the double homogeneous coordinate $(X:Z)\times(Y:W)$; moreover, the curve $X_{\lambda}$ is defined by the equation (\ref{*}) in \S\ref{sec-ASMcurves}.
By Proposition \ref{prop-ramifications1} (4), we may assume that the coordinates are chosen in such a way that, for any $u\in\F_q$,
$$
\begin{array}{ccl}
\varphi^{-1}((u:1)\times(1:0))&=&\{\gamma(u)\,|\,\gamma\in\Gamma\},\\
\varphi^{-1}((1:0)\times(u:1))&=&\{\gamma(t/u)\,|\,\gamma\in\Gamma\}.
\end{array}
$$
Notice that the elements of the right-hand sets are precisely parametrized by $\Gamma$, that is, for example, $\gamma(u)=\gamma'(u)$ for some $u\in\F_q$ implies $\gamma=\gamma'$, since $u$ cannot be a fixed point of non-trivial elements of $\Gamma$ (Remark \ref{rem-NBordinary}).

We view $x=X/Z$ and $y=Y/W$ as functions on $\Omega$ by the map $\varphi$.
\end{se}
\begin{prop}\label{prop-zeropole}\begin{enumerate}
\item For any $v\in\F_q$, the function $x-v=x(z)-v$ vanishes exactly at $\{\gamma(v)\,|\,\gamma\in\Gamma\}$, and each zero is of order $q$.
It has poles exactly at $\{\gamma(t/u)\,|\,\gamma\in\Gamma,\ u\in\F_q\}$, and they are simple poles.

\item For any $v\in\F_q$, the function $y-v=y(z)-v$ vanishes exactly at $\{\gamma(t/v)\,|\,\gamma\in\Gamma\}$, and each zero is of order $q$.
It has poles exactly at $\{\gamma(u)\,|\,\gamma\in\Gamma,\ u\in\F_q\}$, and they are simple poles.

\end{enumerate}

\end{prop}

\begin{proof}
Clearly, $(0:1)\times(1:0)$ is the only point on $X_{\lambda}$ with $x=0$.
Hence we have $x^{-1}(0)=\{\gamma(0)\,|\,\gamma\in\Gamma\}$.
Each zero of $x$ is of order $q$, since the covering map $\phi\colon X_{\lambda}\rightarrow X_0$ ramifies completely at $(0:1)\times(1:0)$.
The proof of the other assertions is similar.
\end{proof}

\begin{thm} \label{mainthm} Set $Q:=q-1$. 
 The algebraic deformation parameter $\lambda$ can be written as an infinite product in the analytic deformation parameter $t$ as follows: 
$$
\lambda(t)=t^{Q^2+1}(1-t^Q)^2\prod_{\gamma\in \Gamma-\{1\}}\bigg[-\frac{(1-\gamma(\infty))(1-\gamma(t)^{Q})}{1-\gamma(0)^q}\frac{(t-\gamma(0))(t^{Q}-\gamma(1)^{Q})}{t^q-\gamma(\infty)^q}\frac{\gamma(\infty)^{Q}\gamma(t)^{Q^2}}{\gamma(1)^{qQ}}\bigg].
$$

\end{thm}

\begin{proof}
The idea of the proof is to write the functions $x(z)-v$ and $y(z)-v$ (for $v \in \F_q$) as infinite products. For example, we want to write $$
\mbox{``\ } x(z)-v=\kappa_v\prod_{\gamma\in\Gamma}\prod_{u\in\F_q}\frac{z-\gamma(v)}{z-\gamma(t/u)} \mbox{\  ''}$$
for some non-zero constant $\kappa_v\in K^{\times}$, matching up zeros and poles of the left hand side with those of the right hand side in $\Omega$, see Proposition \ref{prop-zeropole}. 
The right hand side doesn't really make sense if $\gamma=1$ and $u=1$, but we can replace it by 
$$
x(z)-v=\kappa_v\frac{z^q-v}{z^{q-1}-t^{q-1}}\prod_{\gamma\neq 1}\bigg[\frac{z-\gamma(v)}{z-\gamma(\infty)}\prod_{u\neq 0}\frac{z-\gamma(v)}{z-\gamma(tu)}\bigg]
$$
for some non-zero constant $\kappa_v\in K^{\times}$.

Note that the factor $$\frac{z^q-v}{z^{q-1}-t^{q-1}}$$ in the above product is indeed a function with a zero at $z=v$ of order $q$ and single poles at $z=t/u$ ($u \neq 0$) and $z=\infty$ (since the function has degree $q-(q-1)=1$ as a rational function of $z$). 

The right hand side is a uniformely convergent expression in $z$ on compact subsets of the ordinary set, cf.\ \cite{GvdP}, pp.\ 44-47. Since it has the same zeros and poles as $x(z)-v$, it differs from this function by a constant $\kappa_v$ (since the quotient curve is compact, so has no globally holomorphic functions).  
One can determine this $\kappa_v$, since we know $x(u)=u$ for any $u\in\F_q$.
We calculate $\kappa_v$ by computing $x(0)$ for $v\neq 0$ and $x(1)$ for $v=0$. We may indeed substitute these values for $z$ since they are in the ordinary set, cf.\ Remark \ref{rem-NBordinary}. We find:
$$
\kappa_v=\left\{\begin{array}{ll}
{\displaystyle -t^{q-1}\prod_{\gamma\neq 1}\bigg[\frac{\gamma(\infty)}{\gamma(v)}\prod_{u\neq 0}\frac{\gamma(tu)}{\gamma(v)}\bigg]}&(v\neq 0)\\
{\displaystyle (1-t^{q-1})\prod_{\gamma\neq 1}\bigg[\frac{1-\gamma(\infty)}{1-\gamma(0)}\prod_{u\neq 0}\frac{1-\gamma(tu)}{1-\gamma(0)}\bigg]}&(v=0).
\end{array}\right.
$$

We do a similar computation with the $y$-functions: 
$$
\mbox{`` \ }y(z)-v=\kappa'_v\prod_{\gamma\in\Gamma}\prod_{u\in\F_q}\frac{z-\gamma(t/v)}{z-\gamma(u)}, \mbox{'' \ }
$$
that is, for $v\neq 0$, 
$$
y(z)-\frac{1}{v}=\kappa'_{v^{-1}}\frac{z^q-t^qv}{z^q-z}\prod_{\gamma\neq 1}\bigg[\frac{z-\gamma(tv)}{z-\gamma(0)}\prod_{u\neq 0}\frac{z-\gamma(tv)}{z-\gamma(u)}\bigg],
$$
and
$$
y(z)=\kappa'_0\frac{1}{z^q-z}\prod_{\gamma\neq 1}\bigg[\frac{z-\gamma(\infty)}{z-\gamma(0)}\prod_{u\neq 0}\frac{z-\gamma(\infty)}{z-\gamma(u)}\bigg].
$$
Now the constants $\kappa'_v$ are calculated by evaluating $y(\infty)$ $(v\neq 0)$ and $y(t)$ $(v=0)$:
$$
\kappa'_v=\left\{\begin{array}{ll}
{\displaystyle -v}&(v\neq 0)\\
{\displaystyle (t^q-t)\prod_{\gamma\neq 1}\bigg[\frac{t-\gamma(0)}{t-\gamma(\infty)}\prod_{u\neq 0}\frac{t-\gamma(u)}{t-\gamma(\infty)}\bigg]}&(v=0).
\end{array}\right.
$$

In principle, one can now find an infinite product for $\lambda$ by multiplying together these functions. But rather, we first do a further simplication of the infinite products by using the identities in \S\ref{sec-ASMcurves}. They tell us that the elements $\mu_u$ (for $u\in\F^{\times}_q$) normalize $\Gamma$ (cf.\ Remark \ref{rem-normalizer}), that is, we have $\Gamma\mu_u=\mu_u\Gamma$.
By this, for example, we have
$$
\prod_{\gamma\neq 1}\gamma(tu)=\prod_{\gamma\neq 1}u\gamma(t).
$$
By using this fact, one can now simplify the formula as follows:
$$
\kappa_v=\left\{\begin{array}{ll}
{\displaystyle -t^{q-1}\prod_{\gamma\neq 1}\frac{\gamma(\infty)\gamma(t)^{q-1}}{v\gamma(1)^q}}&(v\neq 0)\\
{\displaystyle (1-t^{q-1})\prod_{\gamma\neq 1}\frac{(1-\gamma(\infty))(1-\gamma(t)^{q-1})}{1-\gamma(0)^q}}&(v=0),
\end{array}\right.
$$
and
$$
\kappa'_0=(t^q-t)\prod_{\gamma\neq 1}\frac{(t-\gamma(0))(t^{q-1}-\gamma(1)^{q-1})}{t^q-\gamma(\infty)^q}.
$$

We are interested in evaluating the functions $x^q-x$ and $y^q-y$.
By the above formula, we now see that $$\lambda=(x^q-x)(y^q-y)=\prod_{v\in\F_q}(x-v)(y-v)=\prod_{v\in\F_q}\kappa_v\kappa'_v$$ is given by
\[
\lambda(t)=t^{(q-1)^2+1}(1-t^{q-1})^2 p(t),
\]
where
\begin{eqnarray} \label{ldef}
\nonumber 
p(t) &= & 
\prod_{\gamma\neq 1}\bigg[-\frac{(1-\gamma(\infty))(1-\gamma(t)^{q-1})}{1-\gamma(0)^q}\frac{(t-\gamma(0))(t^{q-1}-\gamma(1)^{q-1})}{t^q-\gamma(\infty)^q}\frac{\gamma(\infty)^{q-1}\gamma(t)^{(q-1)^2}}{\gamma(1)^{q(q-1)}}\bigg]
\end{eqnarray}

Observe that the auxiliary variable $z \in \Omega$ has disappeared from the final result. This proves the product expansion of $\lambda$ as a function of $t$. 

\end{proof}

\section{ Analytic behaviour of $\lambda(t)$}

In this section, we study the analytic properties of $\lambda$ as a function of $t$. It is important to realise that we have, up to now, considered infinite product `theta functions' in the variable $z \in \Omega$, where $t$ occurs as a `parameter' in both the group elements over which the infinite product runs, and as argument for the corresponding M\"obius transformations. We now have to switch viewpoint and consider the functional dependence on $t$. However, the convergence properties of the used theta series products in $z$ imply that the infinite product on the right hand side in the main formula from Theorem \ref{mainthm} is absolutely convergent \emph{for a fixed value of $t$ with $0<\abs{t}<1$}, so this expression defines 
$$ \lambda \, : \, \Delta_0 \rightarrow \Delta_0 \, : \, t \mapsto \lambda(t) $$
as a \emph{pointwise continuous} function from the punctured open unit disk to itself. But there is more:

\begin{thm} \label{analytic}
The function $\lambda(t)$ is a holomorphic function from $\Delta_0 \rightarrow \Delta_0$, and it extends across zero to a holomorphic function from the open disk $\Delta$ to itself.
\end{thm}
\begin{proof}
For every specific value of $t=t_0$ with $0<|t_0|<1$, we know that the infinite product defining $\lambda(t_0)$ is absolutely convergent. We define rational functions $g_n, n>0$ by
\begin{equation} \label{gnt}
g_n(t)=\prod_{\ell(\gamma)=n} \bigg[-\frac{(1-\gamma(\infty))(1-\gamma(t)^{q-1})}{1-\gamma(0)^q}\frac{(t-\gamma(0))(t^{q-1}-\gamma(1)^{q-1})}{t^q-\gamma(\infty)^q}\frac{\gamma(\infty)^{q-1}\gamma(t)^{(q-1)^2}}{\gamma(1)^{q(q-1)}}\bigg],
\end{equation}
where $\ell$ is the word length function on the free group $\Gamma$ for the given generators $\delta_{u,v}$. 

It suffices to prove that the product $\prod g_n(t)$ is uniformly convergent on closed annuli in $\Delta_0$. So fix constants $0< k_1,k_2 < 1$ and assume $k_1 \leq |t| \leq k_2$.  

The rational functions $g_n(t)$ have coefficients in ${\bf F}_q$, in particular, of absolute value 1 (if they are not zero). Therefore, if $|t|<1$, we have $|g_n(t)-1|=|t|^{e_n}$ for some integer $e_n$ that only depends on $n$ and not on $t$. (As a matter of fact, $e_n$ is the difference between the degrees of the lowest non-vanishing terms in numerator and denominator of $g_n(t)-1$.) The absolute convergence at a particular point $t_0$ with $|t_0|<1$ implies $|g_n(t_0)-1|=|t_0|^{e_n} \rightarrow 0$ as $n \rightarrow + \infty$, and this gives that $e_n \rightarrow +\infty$ as $n \rightarrow +\infty$. Hence we find that for all $k_1 \leq |t| \leq k_2$, $|g_n(t)-1|<c_n:=k_2^{e_n}$, with $c_n \rightarrow 0$ \emph{uniformly} in $t$. Therefore, the infinite product $g_n(t)$ is also uniformly convergent to a holomorphic function on $\Delta_0$. 
As $\lambda=\lambda(t)$ for $t \in \Delta_0$ corresponds to a Mumford curve, $\lambda(t)$ takes values in $\Delta_0$, and in particular it is bounded as a function of $t$. Hence it extends across zero to a holomorphic map on $\Delta_0$ by Lemma \ref{analyticlem}, (3). 
\end{proof}

\begin{rem}
For a general family of Schottky groups that depend algebraically on a deformation parameter $t$, the corresponding algebraic moduli depend rigid-analytically on $t$. In general, one can use the theory of non-archimedean theta functions (cf.\ \cite{ManDrin}) to construct coordinate functions on the algebraic curve, by identifying their zeros and poles. After specialisation of the argument to enough ordinary points, one finds that the algebraic moduli can be given explicitly as algebraic functions of analytic functions of $t$ that are infinite products of rational functions in $t$. This is the general pattern of the proof of our main results above. \end{rem}

\section{$\lambda(t)$ as a function of moduli}

We now study what the function $\lambda(t)$ does if we consider it to depend on the moduli of curves, instead of the parameters $\lambda$ and $t$ itself. For that, we will first determine when two ASM-curves are isomorphic, and when two ASM-groups are conjugate. 

\begin{thm} \label{subraocurves}
Two ASM-curves  $X_{\lambda_i}:(x^q-x)(y^q-y)=\lambda_i$, $i=1,2$ $\lambda_i \in K^*$ are isomorphic if and only if $\lambda_1/\lambda_2 \in \F^*_q$.
For $\zeta \in \F^*_q$, an isomorphism $\psi$ between $X_\lambda$ and $X_{\zeta\lambda}$ is given by $\psi(x)=\zeta x$, $\psi(y)=y$. In particular, since $\F_q \subset K$, isomorphisms of ASM-curves are defined over $K$.  
\end{thm}
\begin{proof}

Let $F_\lambda$ denote the function field of the ASM-curve $$X_\lambda:(x^q-x)(y^q-y)=\lambda.$$
The group $\mathrm{Aut} (X_\lambda)$ of the curve  is generated  by the following elements:
$\tau_{a,b}(x,y)=(x+a,y+b)$ where $a,b \in \F_q$,
$\sigma_1(x,y)=(y,x)$, $\sigma_2(x,y)=(\epsilon x, \epsilon^{-1} y)$, where $\epsilon$ is a primitive $(q-1)$-th 
root of $1$, cf.\ \cite[th.7]{Val-Mad:80}. 

\begin{df} Call a subfield $F$ of $F_\lambda$ \emph{good} if it is a rational subfield of $F_\lambda$ so that $F_\lambda/F$ is an extension with Galois group $\mathrm{Gal}(F_\lambda/F)$ isomorphic to $(\F_q,+)$ as a subgroup of $\mathrm{Aut} (X_\lambda)$. 
\end{df}

\begin{lem}
The fields $K(x)$ and $K(y)$ in $F_\lambda$ are the only good fields $F$ with $G=\Gal(F_\lambda/F)$ satisfying: \begin{enumerate} \item[(a)] $G$ is stable by $\sigma_2$-conjugation, \item[(b)] The action of conjugation by $\sigma_2$ is transitive on $G-\{1\}$ and \item[(c)] there does not exist a subgroup $H$ of $G$ of order $p$ whose $\sigma_1$-conjugate lies in $G$. \end{enumerate}
\end{lem}

\begin{proof}[Proof of Lemma.] Let $F$ be good. Then $G$ is generated by certain $\tau_{a_i,b_i}$ for $i=1,\dots,m$ (for a proof, use that $\Aut(X_\lambda)$ has $E \times E'$ as unique Sylow-$p$ subgroup). We consider three cases:
\begin{enumerate} 
\item \emph{There exists an $i$ such that $a_ib_i \neq 0$.}\\ Set $a=a_i, b=b_i$. Now $a^{-1} b \in \F_q^*$ and since $\epsilon$ generates $\F_q^*$, there exists an integer $k$ such that $a^{-1}b = \epsilon^k$. One computes that $$\sigma_2^k \tau_{a,b} \sigma_2^{-k} = \tau_{\epsilon^k a, \epsilon^{-k} b} = \tau_{b,a} = \sigma_1 \tau_{a,b} \sigma_1.$$ Now if $G$ is $\sigma_2$-conjugation stable, we find that $H=\langle \tau_{a,b} \rangle$ is a group of order $p$ in $G$ with a $\sigma_1$-conjugate inside $G$. Hence in this case, condition (a) and (c) cannot be satisfied simultaneously.\\
\item \emph{Some $a_i=0$ or some $b_i=0$, but not all $a_i=0$ or $b_i=0$.}\\ We verify that condition (b) is not satisfied in this case. Suppose without loss of generality that $(a,0)$ and $(a',b)$ occur, for some $b \neq 0$. Then $\sigma_2 \tau_{a,0} \sigma_2^{-1} = \tau_{\epsilon a, 0} \neq \tau_{a',b}$. \\
\item \emph{All $a_i=0$ or all $b_i=0$.}\\ Now the fixed field of $G$ is $K(x)$ in the first case and $K(y)$ in the second case. We only treat the first case, as the second one is symmetric after interchange of $a_i$ and $b_i$. We verify the indicated properties: for $a=a_i \neq 0$, $\sigma_2 \tau_{a,0} \sigma_2 = \tau_{\varepsilon a,0} \in G$ and since the action of $\langle \varepsilon \rangle = \F_q^*$ is transitive on $\F_q^*$, we find the same for the action of $\sigma_2$ on $G-\{1\}$. For (c), the subgroups $H$ of order $p$ in $G$ are generated by some $\tau_{a,0}$, and $\sigma_1$ maps these to $\tau_{0,a} \notin G$. 
\end{enumerate}

\end{proof}

Now consider two curves  $X_{\lambda_i}:(x^q-x)(y^q-y)=\lambda_i$, $i=1,2$ and $\lambda_i \in K^*$ with
corresponding function fields $F_{\lambda_i}$
 and let $\psi:F_{\lambda_1}\rightarrow F_{\lambda_2}$ 
be an isomorphism. 
The map
$$
\mathrm{Aut}(F_{\lambda_1}) \rightarrow \mathrm{Aut}(F_{\lambda_2}) \ : \ 
\sigma \mapsto  \psi \sigma \psi^{-1},
$$
is an isomorphism of the corresponding  automorphism groups. This implies in particular that the subgroup structure of the automorphism group is preserved by $\psi$. 
Let $A \cong (\F_q,+)$ be the Galois group $\mathrm{Gal}(F_{\lambda_1}/K(x))$. Recall that $K(x)$ is stable under $\sigma_2$. 
The element $\psi(x)$ generates a rational subfield of the function field $F_{\lambda_2}$ of 
the curve $X_{\lambda_2}$ and  $K(\psi(x))= F_{\lambda_2}^{\psi A \psi^{-1}}$, so it is good in $F_{\lambda_2}$. Since $\psi$ preserves the structure of the automorphism group, $K(\psi(x))$ is good and inherits properties (a), (b) and (c) from the lemma above. By that lemma,  
$K(\psi(x))=K(x)$ or $K(\psi(x))=K(y)$. By composing, if necessary, with $\sigma_1$ to interchange $x$ and $y$, we can asssume that $K(\psi(x))=K(x)$. Therefore 
$\psi(x)=g(x)$ for some invertible fractional transformation $g \in \PGL(2,K)$. 
Since  $\psi$ is a morphism of fields, it preserves the algebraic relation defining  $X_{\lambda_1}$. Hence if we set $y'=\psi(y)$, we find
\[
 \left({y'}^q-y' \right)(g(x)^q-g(x))=\lambda_1
\]
Thus $y'$ is a generator of the Artin-Schreier extension $F_{\lambda_2}/K(g(x))$
and according to Hasse \cite[eq. 3']{Hasse34} and Pries \cite[Lemma 2.4]{Pries} it is related to the generator $y$
by
\begin{equation} \label{ar1}
y'=\zeta y +d, \qquad \frac{\lambda_2}{x^q-x}=   \zeta \cdot \frac{\lambda_1}{g(x)^q-g(x)} + d^q-d
\end{equation}
for some $\zeta \in \F_q^*$ and $d \in K(x)$. 
We find in particular that
\begin{equation} \label{blub}  \frac{\lambda_2}{x^q-x} -  \zeta \cdot \frac{\lambda_1}{g(x)^q-g(x)} = d^q-d, \end{equation}  for some $d \in K(x)$. Now if the function $d^q-d$ has any poles, those are of order divisible by $q$, whereas the left hand side of the expression has at most simple poles at  $x \in \F_q$ and $g^{-1}(\F_q)$. We conclude that $d$ doesn't have any poles, so  the left hand side of (\ref{blub}) doesn't have any poles either. But the two individual terms have poles with respective constant residue $\lambda_2$ and $\zeta \lambda_1$, so in corresponding poles, the residues have to cancel. Hence the desired result $\zeta \lambda_1 = \lambda_2$. 

Conversely, if $\lambda_1/\lambda_2=\zeta \in \F_p^*$ then the transformation
 $\psi(x)=\zeta x$, $\psi(y)=y$ makes the function fields $F_{\lambda_i}$, $i=1,2$ isomorphic.
\end{proof}

On the analytic side, we independently have the following result:

\begin{thm} \label{asmconj}
Two Artin-Schreier-Mumford curves $X_{\lambda(t_i)}$, $i=1,2$ uniformized by
 Schottky groups $\Gamma(t_i)$ for $t_i \in \Delta_0$ are equivariantly isomorphic if and only if the groups $\Gamma(t_i) \ (i=1,2)$ are conjugate subgroups of $\PGL(2,K)$, if and only if $t_1=\zeta t_2$ for an 
 element $\zeta \in \F_q^*$.
\end{thm}
\begin{proof}
In this proof, we will refer to the features of an individual ASM-curve by the letters used up to now ($N,\Gamma,\epsilon_u$, etc.),  but if we refer to the dependence of such features on an analytic deformation parameter $t$, we will indicate this by adding $t$ in brackets (so $N(t), \Gamma(t), \epsilon_u(t)$ etc.) Note that in this proof, the notation $g(t)$ for a matrix $g$ does \emph{not} mean the evaluation of the fractional linear transformation $g$ at the point $t$. 

Assume that the Artin-Schreier-Mumford curves $X_{\lambda(t_i)}$ are isomorphic. 
According to \cite[Corolary 4.11]{Mumford:72}
this is equivalent to the fact that the groups $\Gamma(t_i)$ are conjugate in $\PGL(2,K)$, i.e. there is a $\gamma \in \PGL(2,K)$ 
such that $\gamma \Gamma(t_1) \gamma^{-1}=\Gamma(t_2)$. So assume this. Let $N(t_i)$ denote the normalizer of 
$\Gamma(t_i)$ in $\PGL(2,K)$. It follows from the definition of normalizer that $N(t_1)$ and $N(t_2)$ are also conjugate by the same $\gamma$: $\gamma N(t_1) \gamma^{-1} =N(t_2)$. 

For $i=1,2$, we can normalize our ASM-groups so that $N(t_i) \subseteq \PSL(2,\F_q[t_i^{-1}])$, by multiplying $\epsilon'_v$ by $t_i^{-1}$. Now the particular element 
$A_i:=\epsilon_1 \epsilon'_1$ has trace $2+t_i^{-1}$, and trace is conjugacy invariant in $\PSL$. (If one wishes to work in $\PGL$ instead, one can use the conjugacy invariant trace$^2/\det$.) Hence the relations $\gamma a_i \gamma^{-1} \in N(t_{i+1})$ for $i=1,2 \mbox{ mod } 2$,  imply that $2+t_1^{-1} \in \F_q[t_2^{-1}]$ and $2+t_2^{-1} \in \F_q[t_1^{-1}]$. This implies that $t_1^{-1} = \alpha t_2^{-1} + \beta$ for some elements $\alpha \in \F_q^*, \beta \in \F_q$. (We thank Andreas Schweizer for showing us this argument.) 

This already reduces the number of conjugate groups to a finite set. We will now show that $\beta =0$ if the groups are conjugate. 

For this, we first observe that all elements of order $p$ in $N$ are conjugate to $\epsilon_1$ \emph{by a conjugacy in $N$}. Indeed, since dividing modulo a free group $\Gamma$ doesn't kill any torsion elements, an element of order $p$ in $N$ is pulled back from $E \times E'$ in diagram (\ref{dia}). Hence it is an element of order $p$ in the free product $E \ast E'$. But any such element is conjugate \emph{in $E \ast E'$} to an element of $E$ or $E'$ (cf.\ \cite{Serre}, 4.3). Now any element in $E$ (resp.\ $E'$) equals $\epsilon_1$ (resp.\ $\epsilon_1'$) after conjugation with an appropriate $\mu_a \in N$. Finally, it suffices to note that $\epsilon'_1$ is conjugate to $\epsilon_1$ in $N$ via $\tau$ by definition. 

The element $\gamma \epsilon_1 \gamma^{-1} \in N(t_2)$ can hence be conjugated back to $\epsilon_1$ by an element $\delta \in N(t_2)$. Replacing $\gamma$ by $\delta \gamma$, we can henceforth assume that $\gamma$-conjugation fixes $\epsilon_1$. It is easily computed that this implies $\gamma=\bigl( \begin{smallmatrix} 1 & b \\ 0 & 1 \end{smallmatrix} \bigr)$ for some $b \in K$. 

Now back to the dependence on $t_1,t_2$ with $t_i \in \Delta_0$. We compute $\gamma \delta_{1,1}(t_1) \gamma^{-1} = \bigl( \begin{smallmatrix} * & * \\ t_1^{-2} & * \end{smallmatrix} \bigr)  \in N(t_2)$.  Using the relation $t_1^{-1} = \alpha t_2^{-1} + \beta$, the lower left entry of this matrix has absolute value $1/|t_2|^2$. Now it is also an element of $\Gamma(t_2)$, so a word in $\delta_{u,v}(t_2)$.  One proves that the lower left entry of a word $$w:=\delta_{u_1,v_1}(t_2) \dots \delta_{u_r,v_r}(t_2)$$ has absolute value $1/|t_2|^{2r}$. This can be proven by induction on the word length, by calculating only the dominant term (in absolute value) of the product, which will equal (recall $|t_2|<1$)
$$ w=\begin{pmatrix} a t_2^{-2r}+\dots  & b t_2^{1-2r}+\dots \\ c t_2^{-2r}+\dots & d t_2^{1-2r}+\dots \end{pmatrix}  $$ for some non-zero constants $a,b,c,d$ depending on $u_1,v_1,\dots,u_r,v_r$.

Therefore, the conjugate expression has word length one, i.e., it is of the form $\delta_{u,v}(t_2)$ for some $u,v$. But this has lower left entry exactly equal to $uv^2/t^2_2$, which cannot equal $(\alpha t_2^{-1} + \beta)^2$ for non-zero $\beta$. Hence $\beta=0$. 

One now easily checks that $\Gamma(t_1)$ and $\Gamma(t_2)$ are indeed conjugate if $t_1  = \zeta t_2$ for some $\zeta \in \F_q^*$. Actually, we have $\Gamma(t_1)=\Gamma(\zeta t_2)$ in this case, since $\epsilon'_u(t_1)=\zeta \cdot \epsilon'_{\zeta u}(t_2)$. 

\end{proof}

Consider the following purely algebraic problem:

\begin{prob}
Determine whether or not two subgroups $H_1$ and $H_2$ of a given group $G$ are conjugate in that group.
\end{prob}

In general, this is believed to be a computationally hard problem. For subgroups of linear groups, an almost efficient algorithm is known, cf.\ e.g.\ Roney-Dougal \cite{Roney}. We consider the case where $K$ is a field containing $\F_q$, and $H_1$ and $H_2$ are groups of the form $N_1(t)$ for some $t \in K$, where $N_1(t)$ is the subgroup of $\PGL(2,K)$ generated $\varepsilon_u$ and $\varepsilon'_u$, with 
 $$
\varepsilon_u:=\bigg(\begin{array}{cc}1&u\\ 0&1\end{array}\bigg),\quad\tau:=\bigg(\begin{array}{cc}0&t\\ 1&0\end{array}\bigg),\quad \varepsilon'_u:=\tau\varepsilon_u\tau
$$
for all $u\in\F_q$.

If we consider the case where $K$ is non-archimedean valued, then as an amusing corollary of the above theorem (or rather, of a substatement in its proof), we get the following `analytic' solution to a part of the group theoretical conjugacy subgroup problem for $N(t_1)$ and $N(t_2)$ in $\PGL(2,K)$:
\begin{cor} Let  For $t_1,t_2 \in K$ with $0<|t_i|<1 \ (i=1,2)$, we have that $N_1(t_1)$ and $N_1(t_2)$ are conjugate subgroups of $\PGL(2,K)$ if and only if $t_1/t_2 \in \F_q^*$.

\end{cor}

\begin{proof}
It suffices to observe that conjugacy of $N_1(t_1)$ and $N_1(t_2)$ again implies conjugacy of their normalizers $N(t_1)$ and $N(t_2)$, and use the previous argument.

\end{proof}

We don't know to what extent the statement of the corollary is true for arbitrary $t_1,t_2$ and an arbitrary field $K$, and whether it can be proven in a purely algebraic way.

\begin{rem}
From the previous theorem, we also expect that $\lambda$ should depend, up to an element of $\F_q^*$, only on $t^{q-1}$. One may indeed verify directly from the infinite product that this is the case. Since it is not entirely trivial, we prove it here:
\begin{prop}
For every $|t|<1$ and $\zeta\in \mathbf{F}_q^*$ we have 
\[
\lambda( \zeta t)=\zeta \lambda(t).
\] 
\end{prop}

\begin{proof} 

Since the behaviour of the factors in the product for this transformation is somewhat different, we decompose the function $p(t)$ as $p(t)=\prod\limits_{n=1}^\infty p^{(n)}_1(t) \cdot p^{(n)}_2(t) \cdot p^{(n)}_3(t),$
where 
\begin{equation} \label{p1}
p_1^{(n)}(t):= \prod_{\ell(\gamma)=n}
 \left[-\frac{(1-\gamma(\infty))\gamma(\infty)^{q-1}(t^{q-1}-\gamma(1)^{q-1})}{(1-\gamma(0))\gamma(1)^{q(q-1)}} \right],
\end{equation}
\begin{equation} \label{p2}
p_2^{(n)}(t):= \prod_{\ell(\gamma)=n}
(1-\gamma(t)^{q-1})\gamma(t)^{(q-1)^2},
\end{equation}
\begin{equation} \label{p3}
p_3^{(n)}(t):= \prod_{\ell(\gamma)=n}
\frac{t-\gamma(0)}{t^q-\gamma(\infty)}, 
\end{equation} for $\ell$ the word length function for the generators $\delta_{u,v}$ of $\Gamma$. 
We will verify that each of the factors is invariant under $t \mapsto \zeta t$ for $\zeta \in \F_q^*$. 

A word in the group $\Gamma(t)$ of lenght $n$ in the generators  $
\delta_{u,v,t}$ (where we now indicate the dependence on $t$ by a subscript since we will later evaluate the matrices as fractional transformations at $t$) will be written as $w_{\bar{u}, \bar{v}, t} = \delta_{u_1,v_1,t} \dots \delta_{u_n,v_n,t}, $ where $\bar{u}$ and $\bar{v}$ are two 
vectors $\bar{v}=(v_1,\ldots,v_n)$, $\bar{u}=(u_1,\ldots,u_n)$ in $(\mathbf{F}_q)^n$ with no $u_i=v_i=0$. 
Observe that $\delta_{u,v,\zeta t}=\delta_{u,\zeta^{-1} v,t}$, 
and therefore we also have
$w_{\bar{u},\bar{v},\zeta t}= w_{\bar{u},\zeta^{-1} \bar{v},t}$. 

We now consider the action of $\zeta$ on the three types of factors $p_i^{(n)}, \ i=1,2,3$ separately.

\begin{enumerate} 
\item In eq. (\ref{p1}) we consider the product  over all elements of lenght $n$, i.e., over all 
vectors $\bar{u},\bar{v} \in {\F}_q^n$ with no $u_i=v_i=0$.
Changing $p^{(n)}_1(t)$ to $p^{(n)}_1(\zeta t)$ is equivalent to the action of $\zeta$ on the vector $\bar{v}$, which only permutes the factors, so we find $p^{(n)}_1(\zeta t) =p^{(n)}_1(t)$.

\item We will prove that  $\gamma(\zeta t)=\zeta \gamma(t)$ for any element $\gamma \in \Gamma$. It then follows by a simple calculation that 
$p^{(n)}_2(\zeta t)=p^{(n)}_2(t)$.
Observe that for an element $\gamma=\delta_{u,v,t}$ of length one we have the following evaluation of the corresponding fractional linear transformation at the point $t$:
\[
\delta_{u,v,t}(t)=\frac{t^2+uv t+ u^2v^2 -u^2v(v-1)}{t+uv(v-1)}.
\]
But then 
$
\delta_{u,v,\zeta t}(\zeta t)=\zeta \delta_{\frac{u}{\zeta},v,t}(t).
$
If $\gamma$ is a word of length $n$ then $\gamma(t)$ is computed by the 
composition of words of length one and the desired result follows by induction.

\item Consider now the factor $p^{(n)}_3(t)$.
Let $\gamma=w_{\bar{u},\bar{v},t}$ be a word of length $n$.  
We compute
\[
\frac{\zeta t-w_{\bar{u},\bar{v},\zeta t}(0)}
{(\zeta t)^q-w_{\bar{u},\bar{v},\zeta t} (\infty)}=
\frac{t-w_{\zeta^{-1} \bar{u},\bar{v},t} (0)}{t^q-w_{\zeta^{-1}\bar{u},\bar{v},t}(\infty)},
\]
therefore $p^{(n)}_3(\zeta t)=p^{(n)}_3(t)$.
\end{enumerate}
The desired result follows by combining all the above computations, since $$\lambda(t)= t \cdot t^{(q-1)^2} (1-t^{q-1})^2p(t).$$
\end{proof}

\end{rem}
\begin{thm} \label{mainthm2}
 The function $\Lambda:\, \Delta_0 \rightarrow \Delta_0, :\, T \mapsto \lambda(T^{\frac{1}{q-1}})^{q-1}$ extends to a rigid analytic automorphism of the open unit disk $\Delta$. Every value in $\Delta_0$ corresponds to a unique curve up to isomorphism in the family $X_\lambda$.
\end{thm}

\begin{proof}
Since $\lambda$ is actually a function of $t^{q-1}$ up scaling by an element of $\F_q^*$, we find that $\Lambda$ is a well-defined rigid-analytic map $\Delta_0 \rightarrow \Delta_0$. Since it is bounded, it extends across $0$. Also, since isomorphism of curves is given by the compatible $\F_q^*$-scaling action on $t$ and $\lambda$, the resulting map is one-to-one. 
\end{proof}

\begin{rem}
We find that $\Lambda(T)=T^{(q-1)^2+1}P(T)$ for some function $P$ with $P(0) \neq 0$. Since $\Lambda$ is also one-to-one, this implies that the continuation of $P(T)$, and hence $p(t)$, to $\Delta$ acquires a finite order pole at $t=0$. Since $p(t)$ is a uniformly convergent product of rational functions in $t$ (away from $t=0$), this implies that some denominator in the product expansion has zeros in $t$. We don't know in which factor(s) this phenomenon occurs.                                                                                                                                                                                                                                                                                                                                                                                                                                                                                                                                                                                                                                                                                                                                                                                                                                                                                                                                                                                                                                                                                                                                                                                                                                                                                                                                                                                                                                                                                                                                                                                                                                                                                                                                                                                                                                                                                                                                                                                                                                                                                                                                                                                                                                                                                                     
\end{rem}

%\clearpage
\begin{small}

\end{small}

\end{document}